\documentclass[11pt,a4paper]{amsart}
\usepackage{amsfonts,amsmath,amssymb,indentfirst, amscd}
\usepackage[utf8]{inputenc}
\usepackage{graphicx}

\newcommand{\er}{\mathbb{R}}
\newcommand{\R}{\mathbb{R}}

\newcommand{\e}{\epsilon}

\newtheorem{Theo}{Theorem}[section]
\newtheorem{Rem}[Theo]{Remark}
\newtheorem{Lem}[Theo]{Lemma}
\newtheorem{Prop}[Theo]{Proposition}

\title{On a nonlinear Schr\"odinger system arising in quadratic media}

\author{Ad\'an J. Corcho}
\address{
\newline {\bf Ad\'an J. Corcho}
\newline
Instituto de Matem\'atica
\newline
Universidade Federal do Rio de Janeiro-UFRJ
\newline
Ilha do Fund\~ao, 21941-909. Rio de Janeiro-RJ, Brazil
\newline
Rio de Janeiro-RJ, Brazil.}
\email{adan@im.ufrj.br}

\author{Sim\~ao Correia}
\address{\newline
{\bf Sim\~ao Correia}
\newline
CMAF-CIO and FCUL
\newline
Campo Grande, Edifício C6, Piso 2, 1749-016 Lisboa, Portugal}
\email{sfcorreia@fc.ul.pt}

\author{Filipe Oliveira}
\address{
\newline
{\bf Filipe Oliveira}
\newline
Mathematics Department and CEMAPRE
\newline
ISEG, Universidade de Lisboa
\newline
 Rua do Quelhas 6, 1200-781 Lisboa, Portugal}
\email{foliveira@iseg.ulisboa.pt}

\author{Jorge D. Silva}
\address{\newline
{\bf Jorge Drumond Silva}
\newline
Center for Mathematical Analysis, Geometry and Dynamical Systems,
\newline
Department of Mathematics,
\newline
Instituto Superior T\'ecnico, Universidade de Lisboa
\newline
Av. Rovisco Pais, 1049-001 Lisboa, Portugal.}
\email{jsilva@math.tecnico.ulisboa.pt}

\begin{document}
\maketitle
\begin{abstract}
We consider the quadratic Schr\"odinger system
\begin{displaymath}
\left\{\begin{array}{llll}
iu_t+\Delta_{\gamma_1}u+\overline{u}v=0\\
2iv_t+\Delta_{\gamma_2}v-\beta v+\frac 12 u^2=0,\quad t\in\er,\,x\in \er^d\times \er,
\end{array}\right.
\end{displaymath}
in dimensions $1\leq d\leq 4$ and for $\gamma_1,\gamma_2>0$, the so-called elliptic-elliptic case. We show the formation of singularities and blow-up in the $L^2$-(super)critical case. Furthermore, we derive several stability results concerning the ground state solutions of this system.


\end{abstract}
\section{Introduction}
\noindent
In this paper we consider the quadratic Schr\"odinger system
\begin{equation}
\label{quadratic}
\left\{\begin{array}{llll}
iu_t+\Delta_{\gamma_1}u+\overline{u}v=0\\
2iv_t+\Delta_{\gamma_2}v-\beta v+\frac 12 u^2=0,\quad t\in\er,\,x\in \er^d\times \er,
\end{array}\right.
\end{equation}
where $d\leq 4$, $\beta,\gamma_1,\gamma_2\in\er$ and $\Delta_{\gamma}=\partial_{x_1}^2+\dots+\partial_{x_d}^2+\gamma \partial_{x_{d+1}}^2$.

\medskip

\noindent
This system arises as a model for the interaction of waves propagating in $\chi^{(2)}$ dispersive media. In the case of electromagnetic waves, these media are characterized by a polarization vector $\mathcal{P}$ of the form $$\mathcal{P}=\epsilon_0\chi^{(1)}(\omega_0)\mathcal{E}+\chi^{(2)}(\omega_0)\mathcal{E}^2.$$
Here,  $\epsilon_0$ is the vacuum permittivity, $\mathcal{E}$ represents the electric field and $\omega_0$ its angular frequency (see \cite{Saut} for a rigourous derivation of \eqref{quadratic} from the Maxwell-Faraday equation and Amp\`ere's Law). In fact, the quadratic Schr\"odinger system \eqref{quadratic} governs  the dynamics of propagation in $\chi^{(2)}$ media in other physical contexts, namely in nonlinear optics (see for instance \cite{opt1}, \cite{opt2}, \cite{opt3}). Despite this wide range of applications, and contrarely to the modelation of propagation  in $\chi^{(3)}$ centrosymmetric media, which give rise to the Kerr nonlinearity (and hence to Schr\"odinger equations with cubic nonlinearities), very few mathematical results concerning quadratic systems are available in the literature.

\bigskip

\noindent
Very recently, in \cite{Saut}, a rigorous mathematical study of \eqref{quadratic} was undertaken in the $L^2-$subcritical case ($d\leq 2$). After establishing the global well-posedness of \eqref{quadratic} in $H^1(\er^{d+1})$, the authors turn their attention to localized solutions, deriving conditions for their existence (or non-existence). Furthermore, when $\gamma_1,\gamma_2>0$, the existence of ground states is shown using the concentration-compactness principle due to P.L. Lions \cite{PLL}. Finally, for $d=2$ and $\beta=0$, the authors prove the orbital stability of these ground states.

\bigskip

\noindent
Before stating the main results of the present paper, we continue this introduction by making some considerations about system \eqref{quadratic} and its localized solutions. Using standard methods, for $d\leq 4$ (that is, in the $H^1-$subcritical case) the following local existence result can be obtained:

\begin{Theo}[Local Well-posedness]
\label{Theo:LWP}
Let $d\leq 4$. The IVP \eqref{quadratic} with initial data $(u_0,v_0)\in H^1(\er^{d+1})\times H^1(\er^{d+1})$ admits a unique maximal solution
$$(u,v)\in C\big([0,T^*);H^1(\er^{d+1})\times H^1(\er^{d+1})\big).$$
If $T^*<+\infty$ then
$$\lim_{t\to T^*} \|\nabla u(t)\|_2^2+\|\nabla v(t)\|_2^2=+\infty.$$
\end{Theo}

\bigskip
\noindent
Also, the following quantities are formally conserved by the flow of (\ref{quadratic}):\\
the mass
\begin{equation}
\label{mass}
M(u(t),v(t))=\int_{\er^{d+1}}\Big( |u(t,x)|^2+4|v(t,x)|^2\Big)dx
\end{equation}
and the energy
\begin{multline}
\label{energy}
E(u(t),v(t))=\frac 12\int_{\er^{d+1}} \Big(|\nabla u(t,x)|_{\gamma_1}^2+|\nabla v(t,x)|^2_{\gamma_2}\\
+\beta|v(t,x)|^2-Re(\overline{u}^2(t,x)v(t,x))\Big)dx,
\end{multline}
where 
$$|\nabla f|_{\gamma}^2=|\partial_{x_1} f|_2^2+\dots+|\partial_{x_d} f|_2^2+\gamma |\partial_{x_{d+1}} f|_2^2.$$

\bigskip
\noindent
\begin{Rem}\label{remarkglobal}
For $d=3$ consider the case $\gamma_1=\gamma_2=1$ and let $C_{GN}$ the best constant for the vector-valued Gagliardo-Nirenberg inequality 
$$Re \int\overline{u}^2v \le \|u\|_{L^3}^2\|v\|_{L^3}\le C\left(\int |u|^2 + 4|v|^2\right)^{\frac{1}{2}}\int \Big(|\nabla u|^2 + |\nabla v|^2\Big).$$ 
Then, from \eqref{energy} and \eqref{mass} we have
\begin{equation*}\begin{split}
\int \Big(|\nabla u|^2+&|\nabla v|^2\Big)= 2E(u_0, v_0) -\beta\int |v|^2+ Re\int \overline{u}^2v\\
&\le2E(u_0, v_0)-\beta\int |v|^2 + C_{GN}M^{\frac{1}{2}}(u_0, v_0)\int \Big(|\nabla u|^2 + |\nabla v|^2\Big)\\
&=2E(u_0, v_0) -\frac{\beta}{4}M(u_0, v_0)  + \frac{\beta}{4}\int |u|^2 \\
&\hspace{4cm} + C_{GN}M^{\frac{1}{2}}(u_0, v_0)\int \Big(|\nabla u|^2 + |\nabla v|^2\Big).
\end{split}\end{equation*}
Thus, the local solutions given in Theorem \ref{Theo:LWP} can be extended to any time interval $[0, T]$  for all data verifying $M(u_0,v_0)< 1/C^2_{GN}$. In particular, this condition implies that
\begin{itemize}
\item [(i)] $E(u_0, v_0) >0$ when $\beta \ge 0$,

\medskip 
\item [(ii)] $8E(u_0, v_0) > \beta M(u_0, v_0)$ when $\beta < 0$.
\end{itemize}
Moreover, we can conclude from Lemma \ref{lemma_GN} of Section \ref{S4} that  $1/C^2_{GN}=M(P_0, Q_0)$, where $(P_0, Q_0)$ is a ground state of the system. 
\end{Rem}

\medskip
\noindent
In Section \ref{S2}, using these invariants, we compute two Virial identities which yield the following blow-up results in the $L^2$-critical and supercritical cases. Our theorems generalize prior results obtained in \cite{Tanaka} in the case $\beta=0$ and $\gamma_1=\gamma_2$ for $d=4$.
\begin{Theo}[Blow-up, $\boldsymbol{d=3}$]\label{Th-blowup-1}
 Consider the IVP for system \eqref{quadratic} with $d=3$, $\gamma_1=\gamma_2:=\gamma >0$ and initial data $(u_0,v_0) \in (H^1(\er^4)\cap  L^2(\er^4, |x|^2dx))^2$. Let $$(u,v)\in C\Big([0, T^*);\,   (H^1(\er^4)\cap  L^2(\er^4, |x|^2dx))^2\Big)$$
 be the corresponding maximal solution. Assume in addition that 
 \begin{equation}\label{blowup-1-ha}
 E(u_0,v_0) < 0\quad and\quad  \beta > 0
 \end{equation}
 or
 \begin{equation}\label{blowup-1-hb}
 8E(u_0,v_0) <\beta M(u_0,v_0) \quad and\quad  \beta \le 0.
 \end{equation}
Then $T^*< \infty$ and $\displaystyle \lim_{t\to T^*} \|\nabla u(t)\|_2^2=+\infty.$
\end{Theo}
\begin{Theo}[Blow-up, $\boldsymbol{d=4}$]\label{Th-blowup-2}
	Consider the IVP for system \eqref{quadratic} with $d=4$, $\gamma_1, \gamma_2 >0$ and initial data $(u_0, v_0) \in (H^1(\er^5)\cap  L^2(\er^5, |x_{\perp}|^2dx))^2$, where $x_{\perp}=(x_1,\dots,x_d)$. Let  $$(u,v)\in C\Big([0, T^*);\,   (H^1(\er^5)\cap  L^2(\er^5, |x_{\perp}|^2dx))^2\Big)$$ 
be the corresponding maximal solution. Assume in addition that
	\begin{equation}\label{blowup-2-ha}
	E(u_0,v_0) < 0\quad and\quad  \beta > 0
	\end{equation}
	or
	\begin{equation}\label{blowup-2-hb}
	8E(u_0,v_0) <\beta M(u_0,v_0) \quad and\quad  \beta \le 0.
	\end{equation}
	Then  $T^*< \infty$. 
\end{Theo}

\begin{Rem}
Data satisfying negative energy is a natural blow-up condition for this model. This is not the case for the conditions given in \eqref{blowup-1-hb} and \eqref{blowup-2-hb}. Notice, however, that initial data satisfying this new hypothesis may be easily built: indeed, let $U\in H^1(\er^{d+1})$, $U>0$. Setting $u_0=v_0=\lambda U$, $\lambda\in\er$, a simple computation shows that, for large $\lambda$, one has $8E(u_0,v_0)<\beta M(u_0,v_0)$.
\end{Rem}

\bigskip
\noindent
For $\gamma_1,\gamma_2,\omega,4\omega+\beta>0$ (see \cite{Saut}), the system (\ref{quadratic}) admits localized solutions of the form
\begin{equation}
u(t,x)=P(x)e^{i\omega t},\qquad v(t,x)=Q(x)e^{2i\omega t}.
\end{equation}
The functions $P$ and $Q$ satisfy the system
\begin{equation}
\label{elliptic}
\left\{\begin{array}{llll}
-\omega P+\Delta_{\gamma_1}P+\overline{P}Q=0\\
\\
-(4\omega+\beta)Q+\Delta_{\gamma_2}Q+\frac 12 P^2=0.
\end{array}\right.
\end{equation}
Let us denote by $B$ the set of all {\it bound states}, that is, the set of all solutions $(P,Q)\in H:=H^1(\er^{d+1})\times H^1(\er^{d+1})$ of the stationnary system \eqref{elliptic}.

\bigskip

\noindent
We will say that a bound state $(P,Q)$ is a {\it ground state} if $(P,Q)$ minimizes the action 
\begin{equation}
\label{action}
S(u,v)=E(u,v)+\omega M(u,v)
\end{equation}
among all bound states. That is, denoting by $G$ the set of all ground states, we have
$$(P_0,Q_0)\in G\Leftrightarrow (P_0,Q_0)\in B\textrm{ and }\forall (P,Q)\in B, S(P_0,Q_0)\leq S(P,Q).$$ 
It is not difficult to see that $(P_0,Q_0)\in G$ if and only if $(P_0,Q_0)$ is a minimizer of the problem
\begin{equation}
\label{minimization}
\inf\{I(u,v)=K(u,v)+\omega M(u,v)\,:(u,v)\in\mathcal{W}_{(P_0,Q_0)}\},
\end{equation}
where
$$K(u,v)=\int_{\er^{d+1}} \Big(|\nabla u(t,x)|_{\gamma_1}^2+|\nabla v(t,x)|_{\gamma_2}^2\\
+\beta|v(t,x)|^2\Big)dx$$
and
$$\mathcal{W}_{(P_0,Q_0)}=\Big\{(u,v)\in H\,:\, J(u,v):=Re\int_{\er^{d+1}}\overline{u}^2v=Re\int_{\er^{d+1}}\overline{P_0}^2Q_0\Big\}.$$

In Section \ref{S3} we will show the following instability results concerning ground states in the $L^2$-critical and supercritical cases:
\begin{Theo}[Strong instability]
\label{Theo:SI}
Let $\gamma_1=\gamma_2>0$, $\beta=0$ and $d=3$. Let $B$ be the set of all bound states of \eqref{elliptic}. Then $B$ is unstable in the following sense: given $(P,Q)\in B$, there exists a sequence $X_{0,k}\to (P,Q)$
in $H$ such that, for all $k$, the solution $X_k$ of \eqref{quadratic} with initial data $X_{0,k}$ blows up in finite time. 
\end{Theo}
\begin{Theo}[Weak instability]
\label{inst}
 Let $\gamma_1,\gamma_2>0$, $d=3$ and $\beta\neq 0$ or $d\geq 4$ and $\beta \in \R$. For $(P,Q)\in G$,
 we consider its orbit
$$\Sigma=\{f(\theta,y)[P,Q]:=(e^{i\theta}P(\cdot+y),e^{2i\theta}Q(\cdot+y))\,:\,\theta\in\er,y\in\er^{d+1}\}.$$
Then $\Sigma$
is weakly unstable by the flow of \eqref{quadratic}, in the following sense: there exists $\epsilon>0$ and a sequence $X_{0,k}\to (P,Q)$ in $H$ such that 
\begin{itemize}
\item The solution $X_k(t)$ to \eqref{quadratic} with initial data $X_{0,k}$ is global and bounded in $H$;
\item For all $k$, $T_{k}^*=\sup\{T\geq 0\,:\,\forall t\in[0,T],\,X_k(t)\in \Sigma_{\epsilon}\}<+\infty$, where $\Sigma_{\epsilon}$ is the $\epsilon$-neighbourhood of $\Sigma$.
\end{itemize}
\end{Theo}

\bigskip
\noindent 
In what concerns stability of ground states, the proof in \cite{Saut} follows the argument of Cazenave and Lions for the stability of ground states of the nonlinear Schr\"odinger equation, by showing that the solutions of the minimization problem 
\begin{equation}\label{min_energia}
\inf\{E(u,v): M(u,v)=M(P,Q),\ (P,Q) \in G \}
\end{equation}
are precisely the ground states of \eqref{quadratic}. In \cite{Saut}, it was proven that: such a minimization problem has a solution; the solution is a bound state, and so it has an action larger or equal than any ground state; the solution is actually a ground state, by proving that it has the same action as any given ground state. The first and third steps only require that system \eqref{quadratic} is $L^2$-subcritical, meaning that $d\le 2$. However, to show the second step, the procedure used therein only works for $d=2$ and $\beta=0$.

Recalling some arguments used in \cite{Simao2}, one may actually skip the second step, as long as the energy does not contain any $L^2$ terms (in the present situation, it means that $\beta=0$). The consequence is a more direct approach, presented in Section \ref{S4}, which is also valid for $d=1$:
\begin{Theo}
\label{Theo:St}
Suppose that $d\le2$ and $\beta=0$. Then the set of ground states $G$ is stable with respect to the flow generated by \eqref{quadratic}, that is, for each $\delta>0$, there exists $\epsilon >0$ such that, if $(u_0,v_0)\in H$ satisfies
$$
\inf_{(P,Q)\in G}\|(u_0,v_0)-(P,Q)\|_{H}<\epsilon,
$$
then the solution $(u,v)$ of \eqref{quadratic} with initial data $(u_0,v_0)$ satisfies
$$
\sup_{t\ge 0} \inf_{(P,Q)\in G}\|(u(t),v(t))-(P,Q)\|_{H}<\delta.
$$
\end{Theo}

\begin{Rem}
For $d=3$, the Virial blow-up result and Remark \ref{remarkglobal} imply that $M(u_0,v_0)=M(P,Q)$, where $(P,Q)$ is any ground-state for the system with $\beta=0$, is the threshold for blow-up behaviour. Moreover, when $\beta=0$, using the pseudo-conformal transformation, one may exhibit a blow-up solution with critical mass (in the lines of \cite{weinstein}). For $\beta\neq 0$, the pseudo-conformal transformation is not available. The existence of blow-up solutions with critical mass remains an interesting open problem.
\end{Rem}

\section{Virial identity and blow-up}\label{S2}
We begin this Section by noticing that the system \eqref{quadratic} can be put in the Hamiltonian form
\begin{equation}
\label{ham}
\frac{\partial X}{\partial t}(t)=JE'(X(t)),
\end{equation}
where $J$ is the skew-adjoint operator $\left[\begin{array}{lll}
-i&0\\
0&-\frac i2
\end{array}\right]$ and $X=(u,v)$.

\medskip

\noindent
Using this fact, we will derive two global Virial type identities for system (\ref{quadratic}). 
Instead of using the standard technique based on several integrations by parts to calculate the second derivative in time for the variance of the solutions, we use an interesting method presented in \cite{GR} that allows to formally unders\-tand the evolution of certain real functional along the trajectories of Hamiltonian systems. In subsection \ref{dual-dynamics} we describe the general idea of this procedure applied to the system (\ref{quadratic}). Finally, we use the Virial identities obtained to establish two results about the formation of singularities for system (\ref{quadratic}) based on classical convexity arguments.

\subsection{Dual dynamics for system \eqref{quadratic}}\label{dual-dynamics}
Consider a real functional $G$, defined on a dense subspace $\mathcal{V}$  of $L^2(\mathbb{R}^{d+1})$, with continuous derivatives in $L^2(\mathbb{R}^{d+1})$.  The goal is to study the evolution of $G$ along the trajectories of the dynamical system defined by equation 
\eqref{ham}. 

Recalling that $X(t) =(u(t), v(t))$, the time derivative of $G$ calculated along $X(t)$ is given by
\begin{equation}\label{Derivative-G}
\begin{split}
\frac{d}{dt}G(X(t))= \langle G'(X(t)),\; \tfrac{\partial X}{\partial t}(t)\rangle&= \langle G'(X(t)),\; JE'(X(t))\rangle\\
&:=P(X(t)).
\end{split}
\end{equation}
On the other hand, given $\tilde{X}_0:=(\tilde{u}_0, \tilde{v}_0)$, consider the initial value problem
\begin{equation}\label{IVP-G}
\frac{\partial}{\partial t}\tilde{X}(t)= J G'(\tilde{X}(t)),\quad \tilde{X}(0)=\tilde{X}_0,
\end{equation}
which we suppose to be locally well-posed. Thus,
\begin{equation}\label{Derivative-E}
\begin{split}
\frac{d}{dt}E(\tilde{X}(t))= \langle E'(\tilde{X}(t)),\; \tfrac{\partial \tilde{X}}{\partial t}(t)\rangle & = \langle E'(\tilde{X}(t)),\; JG'(\tilde{X}(t))\rangle\\
& = -\langle G'(\tilde{X}(t)),\; JE'(\tilde{X}(t))\rangle\\
& = -P(\tilde{X}(t)).
\end{split}
\end{equation}
Therefore, the validation at time $t = 0$ yields
\begin{equation}\label{Formula-P}
P(\tilde{X}_0)=-\frac{d}{dt}E(\tilde{X}(t))\Big |_{t=0},
\end{equation}
which determines the evolution of $G$ along the trajectories of \eqref{ham}.

\medskip

\noindent
In what follows, we write $x_{\perp}:=(x_1, x_2, \dots, x_d)$, $\nabla_{\perp}:=(\partial_{x_1},\dots,\partial_{x_d})$ and 
$\boldsymbol{x}^{\gamma}=(x_{\perp}, \gamma x_{d+1})$. Also we decompose the energy \eqref{energy} in the following way:
\begin{equation}\label{energy-decomposition}
E(u,v)=E_{\gamma_1}(u) + E_{\gamma_2}(v) + E_{\beta}(v)-E_{Re}(u,v),
\end{equation}
where
\begin{align*}
&E_{\gamma_1}(u)=\frac{1}{2}\int_{\er^{d+1}}|\nabla u|_{\gamma_1}^2dx=\frac{1}{2}\int_{\er^{d+1}}\big(|\nabla_{\perp}u|^2 + \gamma_1|\partial_{x_{d+1}}u|^2\big)dx,\\
&E_{\gamma_2}(v)=\frac{1}{2}\int_{\er^{d+1}}|\nabla v|_{\gamma_2}^2dx=\frac{1}{2}\int_{\er^{d+1}}\big(|\nabla_{\perp}v|^2 + \gamma_2|\partial_{x_{d+1}}v|^2\big)dx,\\
&E_{\beta}(v)=\frac{\beta}{2}\int_{\er^{d+1}}|v|^2dx,\\
&E_{Re}(u,v)=\frac{1}{2}Re\int_{\er^{d+1}}\bar{u}^2vdx. 
\end{align*}
Finally, we set $M_0:=M(u(\cdot, 0), v(\cdot, 0))$ and $E_0:=E(u(\cdot, 0), v(\cdot, 0))$.

\subsection{Virial type identities}
In this subsection we prove the following Virial identities:
\begin{Prop}[Virial identity]\label{Virial-I}
Let $d=3, 4$ and $$(u_0,v_0) \in (H^1(\er^{d+1})\cap L^2(\er^{d+1},|x|^2dx))^2.$$ Then, the variance 
$$\mathcal{V}(t)=\mathcal{V}(u(t),v(t)):=\frac{1}{2}\int_{\mathbb{R}^{d+1}}\big |x|^2(\, |u(t)|^2 + 4|v(t)|^2)dx$$
is finite on the maximal time interval $[0, T^*)$ and  $\mathcal{V}\in C^2\big([0, T^*)\big)$.\\
Furthermore, the following identities hold:  
\begin{enumerate}
	\item [\bf{(i)}] $\displaystyle \frac{d\mathcal{V}}{dt}(t)=2Im\int_{\er^{d+1}}\big( \boldsymbol{x}^{\gamma_1}\cdot \nabla u\,\overline{u}dx + 2\boldsymbol{x}^{\gamma_2}\cdot\nabla v\,\overline{v}\big)dx$,\; 
	\item [\bf{(ii)}] If $\gamma_1=\gamma_2:=\gamma$, 
	\begin{equation*}
	\frac{d^2\mathcal{V}}{dt^2}(t)=4\int_{\er^{d+1}}\left(|\nabla u|_{\gamma}^2 + |\nabla v|_{\gamma}^2\right)dx-(d+\gamma)Re\int_{\er^{d+1}}\bar{u}^2vdx.
	\end{equation*}
	\item [\bf{(iii)}] In particular,  for $d=3$ and $\gamma_1=\gamma_2=1$,
	$$\displaystyle \frac{d^2\mathcal{V}}{dt^2}(t)= 8E_0-4\beta\int_{\er^4}|v|^2dx.$$ 
\end{enumerate}
\end{Prop}

\begin{Prop}[Transverse Virial identity]\label{Virial-II}
Let $d=3, 4$ and $$(u_0,v_0) \in (H^1(\er^{d+1})\cap L^2(\er^{d+1},|x_{\perp}|^2dx))^2.$$ Then, the transverse variance 
$$\mathcal{V}_{\perp}(t)=\mathcal{V}_{\perp}(u(t),v(t)):=\frac{1}{2}\int_{\mathbb{R}^{d+1}}\big |x_{\perp}|^2(\, |u(t)|^2 + 4|v(t)|^2)dx$$
is finite on the maximal time interval $[0, T^*)$ and  $\mathcal{V}_{\perp}\in C^2\big([0, T^*)\big)$.\\
Furthermore, the following identities hold:  
	\begin{enumerate}
		\item[\bf{(i)}] $\displaystyle \frac{d\mathcal{V}_{\perp}}{dt}(t)=2Im\int_{\er^{d+1}}\big(x_{\perp}\cdot \nabla_{\perp} u\,\overline{u}+ 2x_{\perp}\cdot\nabla_{\perp} v\,\overline{v}\big)dx$.
		\item[\bf{(ii)}] 
		$\displaystyle \frac{d^2\mathcal{V}_{\perp}}{dt^2}(t)=4\int_{\er^{d+1}}\left(|\nabla_{\perp}u|^2 + |\nabla_{\perp}v|^2\right)dx -d Re\int_{\er^{d+1}}\bar{u}^2vdx$.
		\item[\bf{(iii)}] In particular, for $d=4$, we have
		$$\displaystyle \frac{d^2\mathcal{V}_{\perp}}{dt^2}(t)=8E_0-4\beta\int_{\er^5}|v|^2dx
		-4\int_{\er^5}\left(\gamma_1|\partial_{x_{d+1}}u|^2 + \gamma_2|\partial_{x_{d+1}}v|^2 \right)dx.$$ 
	\end{enumerate}
\end{Prop}

\subsection*{Proof of Proposition \ref{Virial-I}.} 

{\bf Proof of assertion (i):} We formally apply the technique of dual dynamics to the functional $G(u,v):=\mathcal{V}(u,v)$. The corresponding IVP \eqref{IVP-G} for this functional is defined by
\begin{equation*}
\begin{cases}
\tilde{u}_t=-i|x|^2\tilde{u},& \tilde{u}(0)=\tilde{u}_0,\\
\tilde{v}_t=-2i|x|^2\tilde{v},& \tilde{v}(0)=\tilde{v}_0,
\end{cases}
\end{equation*}
whose solution is given by
$$\big(\tilde{u}(x,t), \tilde{v}(x,t)\big)=\big(e^{-i|x|^2t}\tilde{u}_0, e^{-2i|x|^2t}\tilde{v}_0\big).$$
Then, from \eqref{Formula-P}, we get
\begin{equation*}
\begin{split}
P(\tilde{u}_0, \tilde{v}_0)&=-\frac{d}{dt}E(e^{-i|x|^2t}\tilde{u}_0, e^{-2i|x|^2t}\tilde{v}_0)\Big |_{t=0}\\
&=-\frac{d}{dt}E_{\gamma_1}(e^{-i|x|^2t}\tilde{u}_0)\Big |_{t=0} -\frac{d}{dt}E_{\gamma_2}(e^{-2i|x|^2t}\tilde{v}_0)\Big |_{t=0}\\
&=2Im\int_{\er^{d+1}}\boldsymbol{x}^{\gamma_1}\cdot \nabla\tilde{u}_0\,\overline{\tilde{u}_0}dx +4Im\int_{\er^{d+1}}\boldsymbol{x}^{\gamma_2}\cdot \nabla\tilde{u}_0\,\overline{\tilde{u}_0}dx,
\end{split}
\end{equation*}
since $E_{\beta}(\tilde{v})$ and $E_{Re}(\tilde{u},\tilde{v})$ are independent of time. Thus, it follows from \eqref{Derivative-G} that
\begin{equation}\label{variance-derivative}
\frac{d\mathcal{V}}{dt}(t)=2Im\int_{\er^{d+1}}\big( \boldsymbol{x}^{\gamma_1}\cdot \nabla u\,\overline{u}dx + 2\boldsymbol{x}^{\gamma_2}\cdot\nabla v\,\overline{v}\big)dx\\
\end{equation}
as claimed in (i).

\medskip
{\bf Proof of assertion (ii):} To prove (ii), we choose instead
$$G(u, v):=Im\int_{\er^{d+1}}\big( 2\boldsymbol{x}^{\gamma}\cdot \nabla u\,\overline{u}dx + 4\boldsymbol{x}^{\gamma}\cdot\nabla v\,\overline{v}\big)dx.$$
The corresponding  IVP \eqref{IVP-G} is now
\begin{equation*}
\begin{cases}
\tilde{u}_t=-4\boldsymbol{x}^{\gamma}\cdot \nabla \tilde{u} - 2(d+\gamma)\tilde{u},& \tilde{u}(0)=\tilde{u}_0,\\
\tilde{v}_t=-4\boldsymbol{x}^{\gamma}\cdot \nabla \tilde{v} - 2(d+\gamma)\tilde{v},& \tilde{v}(0)=\tilde{v}_0,
\end{cases}
\end{equation*}
so that
\begin{align*}
&\tilde{u}(x,t)=e^{-2(d+\gamma)t}\tilde{u}_0(e^{-4t}x_{\perp}, e^{-4\gamma t}x_{d+1}),\\ \intertext{and}
&\tilde{v}(x,t)=e^{-2(d+\gamma)t}\tilde{v}_0(e^{-4t}x_{\perp}, e^{-4\gamma t}x_{d+1}).
\end{align*}
Now we proceed with the computation of $\displaystyle P(\tilde{u}_0, \tilde{v}_0)=-\frac{d}{dt}E(\tilde{u}, \tilde{v})\Big |_{t=0}$. Using the change of variables  $(x_{\perp}, x_{d+1})=(e^{4t}y_{\perp},\, e^{4\gamma t}y_{d+1})$, we get
\begin{align*}
E(\tilde{u},\tilde{v})=&\frac{1}{2}\int_{\er^{d+1}}\left(e^{-8t}|\nabla _{\perp}\tilde{u}_0(y)|^2 + \gamma e^{-8\gamma t}|\partial_{x_{d+1}}\tilde{u}_0(y)|^2 \right)dy\\
&+\frac{1}{2}\int_{\er^{d+1}}\left(e^{-8t}|\nabla _{\perp}\tilde{v}_0(y)|^2 + \gamma e^{-8\gamma t}|\partial_{x_{d+1}}\tilde{v}_0(y)|^2 \right)dy\\
&+\frac{\beta}{2}\int_{\er^{d+1}}|\tilde{v}_0(y)|^2dy+\frac{e^{-2(d+\gamma)t}}{2}Re\int_{\er^{d+1}}{\bar{\tilde{u}}}_0^2(y)\tilde{v}_0(y)dy.
\end{align*}
Finally, we conclude that
\begin{equation*}\begin{split}
P(\tilde{u}_0,\tilde{v}_0)&=-\frac{d}{dt}E(\tilde{u}, \tilde{v})\Big |_{t=0}\\
&=4\int_{\er^{d+1}}\left(|\nabla _{\perp}\tilde{u}_0|^2 + \gamma^2|\partial_{x_{d+1}}\tilde{u}_0|^2 + |\nabla _{\perp}\tilde{v}_0|^2 + \gamma^2|\partial_{x_{d+1}}\tilde{v}_0|^2\right)dy\\
&\quad - (d+\gamma)Re\int_{\er^{d+1}}{\bar{\tilde{u}}}_0^2\tilde{v}_0dy, 
\end{split}\end{equation*}
which implies (ii). 

\medskip
{\bf Proof of assertion (iii):} The identity is an immediate consequence of (ii) combined with the conservation of the energy \eqref{energy}.\hfill$\blacksquare$

 \subsection*{Proof of Proposition \ref{Virial-II}}
{\bf Proof of assertion (i):} The proof is similar as the one performed for the case (i) in Proposition \ref{Virial-I} and follows without major changes. 

\medskip
{\bf Proof of assertion (ii):} Here we take  $G$ defined by 
$$G(u, v):=Im\int_{\er^{d+1}}\big( 2x_{\perp}\cdot \nabla_{\perp} u\,\overline{u}dx + 4x_{\perp}\cdot\nabla_{\perp} v\,\overline{v}\big)dx.$$
In this case, the corresponding  IVP \eqref{IVP-G} is written as follows:
\begin{equation*}
\begin{cases}
\tilde{u}_t=-4x_{\perp}\cdot \nabla_{\perp}\tilde{u} - 2d \tilde{u},& \tilde{u}(0)=\tilde{u}_0,\\
\tilde{v}_t=-4x_{\perp}\cdot \nabla_{\perp}\tilde{v} - 2d \tilde{v},& \tilde{v}(0)=\tilde{v}_0,
\end{cases}
\end{equation*}
so that 
$$(\tilde{u}, \tilde{v})=\Big(e^{-2dt}\tilde{u}_0(e^{-4t}x_{\perp}, x_{d+1}),\; e^{-2dt}\tilde{v}_0(e^{-4t}x_{\perp}, x_{d+1})\Big).$$
Using the change of variables $(x_{\perp},\, x_{d+1})=(e^{4t}y_{\perp},\, y_{d+1})$, we have
\begin{equation*}
\begin{split}
E(\tilde{u}, \tilde{v})&=\frac{e^{-8t}}{2}\int_{\er^{d+1}}\left(|\nabla_{\perp}\tilde{u}_0(y)|^2 + |\nabla_{\perp}\tilde{v}_0(y)|^2 \right)dy\\
&\quad + \frac{1}{2}\int_{\er^{d+1}}\left(\gamma_1|\partial_{x_{d+1}}\tilde{u}_0(y)|^2 + \gamma_2|\partial_{x_{d+1}}\tilde{v}_0(y)|^2 \right)dy\\
&\quad + \frac{\beta}{2}\int_{\er^{d+1}}|\tilde{v}_0(y)|^2dy -\frac{e^{-2dt}}{2}Re\int_{\er^{d+1}}{\bar{\tilde{u}}}_0^2(y)\tilde{v}_0(y)dy, 
\end{split}
\end{equation*}
hence
\begin{equation*}
\begin{split}
P(\tilde{u}_0,\tilde{v}_0)&=-\frac{d}{dt}E(\tilde{u}, \tilde{v})\Big |_{t=0}\\
&=4\int_{\er^{d+1}}\left(|\nabla_{\perp}\tilde{u}_0(y)|^2 + |\nabla_{\perp}\tilde{v}_0(y)|^2 \right)dy-dRe\int_{\er^{d+1}}{\bar{\tilde{u}}}_0^2(y)\tilde{v}_0(y)dy,
\end{split}
\end{equation*}
which yields (ii). 

\medskip
{\bf Proof of assertion (iii):} Once again, this last assertion is a particular case of (ii) combined with the conservation of the energy \eqref{energy}.\hfill$\blacksquare$

 \subsection{Proof of the blow-up results}

\noindent{\bf Proof of Theorem \ref{Th-blowup-1}:} By rescaling, it is easy to reduce the problem to the case $\gamma=1$, which can be trated by the classical convexity method, similar to the nonlinear Schr\"odinger equation (see for instance \cite{Cazenave}). The blow-up of $\|\nabla u\|_2^2$ follows from the energy conservation law and the blow-up alternative presented in Theorem \ref{Theo:LWP}. Indeed, from \eqref{energy}, Hölder’s inequality and the Sobolev inequality in dimension $n=4$ it follows that
\begin{equation*}
\begin{split}
\|\nabla_{\gamma}u(\cdot, t)\|_{L^2}^2 + \|\nabla_{\gamma}&v(\cdot, t)\|_{L^2}^2 =2E_0 -\beta\|v(\cdot, t)\|_{L^2}^2+Re\int_{\er^4}\bar{u}^2(\cdot, t)v(\cdot, t)\\
&\le 2E_0 + \frac{|\beta|}{4}M_0 + c \|\nabla u(\cdot, t)\|_{L^2(\er^4)}^2\|v(\cdot, t)\|_{L^2(\er^4)}\\
&\le 2E_0 + \frac{|\beta|}{4}M_0 + c \frac{\sqrt{M_0}}{2}\|\nabla u(\cdot, t)\|_{L^2(\er^4)}^2.
\end{split}
\end{equation*}
\hfill$\blacksquare$

\medskip

\noindent
We finish by noticing that the Virial identity {\bf (iii)} in Proposition \ref{Virial-II} and arguments similar to the ones used in the proof of Theorem \ref{Th-blowup-1} allow us to establish Theorem \ref{Th-blowup-2}.  
\begin{Rem}
Notice that dimensions $d=3,4$ are $L^2$-(super)critical and $H^1$-subcritical. In this situation, the local $H^1\times H^1$  existence theory allows to prove the persistence of solutions in $H^s\times H^s$, $s>\frac{d+1}2$, provided that the initial data has $H^s\times H^s$ regularity. In this framework, one can show the blow-up $$\lim_{t\to T^*}\|v(\cdot,t)\|_{\infty}=+\infty.$$
Indeed, for $d=3$ (a similar computation can be produced for $d=4$):
\begin{equation*}
\begin{split}
\|\nabla_{\gamma}u(\cdot, t)\|_{L^2}^2 + \|\nabla_{\gamma}&u(\cdot, t)\|_{L^2}^2 =2E_0 -\beta\|v(\cdot, t)\|_{L^2}^2+Re\int_{\er^4}\bar{u}^2(\cdot, t)v(\cdot, t)\\
&\le 2E_0 + \frac{|\beta|}{4}M_0 + \|u(\cdot, t)\|_{L^2(\er^4)}^2\|v(\cdot, t)\|_{L^{\infty}(\er^4)}\\
&\le 2E_0 + \frac{|\beta|}{4}M_0 + M_0\|v(\cdot, t)\|_{L^{\infty}(\er^4)}. 
\end{split}
\end{equation*}
\end{Rem}
\section{Instability of ground states}\label{S3}
The proof of Theorem \ref{Theo:SI} follows from Theorem \ref{Th-blowup-1} and from the fact that, given a bound state $(P,Q)$, $E(\lambda P,\lambda  Q)<0$ for $\lambda>1$.  

\medskip

\noindent
We now show the weak instability of ground state solutions to (\ref{quadratic}) in the critical ($d=3$) and supercritical ($d\geq 4$) cases.
\noindent
Let $(P,Q)$ a ground state, that is, a solution of the minimization problem \eqref{minimization}. 

\bigskip

\noindent
Noticing that $\int \overline{P}^2Q\in \er^+$ (see \cite{Saut}), we can show that the orbit of every ground state contains an element $(\tilde{P},\tilde{Q})$, with $\tilde{P},\tilde{Q}>0$. More precisely:
\begin{Prop}
\label{positivas} Let $\gamma_1,\gamma_2>0$ and $n\geq 1$.\\
Then, for every solution $(P,Q)\in H$ of the minimization problem \eqref{minimization}:

\medskip

(i) $(|P|,|Q|)$ is a solution of \eqref{minimization}.\\

(ii) $(|P|,|Q|)$ belongs to the orbit
$$\Sigma=\{e^{i\theta}P(\cdot+y),e^{2i\theta}Q(\cdot+y)\,:\,\theta\in\er,y\in\er^{d+1}\}\textrm{ of $(P,Q)$.}$$

\end{Prop}
\noindent
{\bf Proof of (i):}\\
Let $(P,Q)\in H$ a minimizer and take $(\tilde{P},\tilde{Q})=(|P|,|Q|)$.\\
It is straightforward to see that
$I(\tilde{P},\tilde{Q})\leq I(P,Q)$. Furthermore,
$$\tilde{\mu}=Re \int \tilde{P}^2\tilde{Q}=\int \tilde{P}^2\tilde{Q}=\int |\overline{P}|^2|Q|\geq \Big|\int \overline{P}^2Q\Big|=Re\int \overline{P}^2Q=\mu.$$
Now, assume that $\tilde{\mu}>\mu$. For $\displaystyle\lambda=\Big(\frac{\mu}{\tilde{\mu}}\Big)^{\frac 2n}$, we put
$$(P_{\lambda}(\cdot),Q_{\lambda}(\cdot))=(\lambda^{\frac n2}\tilde{P}(\lambda \cdot),\lambda^{\frac n2}\tilde{Q}(\lambda \cdot)).$$
We get $\displaystyle \int P^2_{\lambda}Q_{\lambda}=\mu$ and $I(P_{\lambda},Q_{\lambda})<I(\tilde{P},\tilde{Q})\leq I(P,Q)$, which contradicts the minimality of $(P,Q)$.\hfill$\blacksquare$

\bigskip

\noindent
{\bf Proof of (ii):}\\
Write $(P,Q)=(|P|e^{\theta_1 (x)},|Q|e^{i\theta_2(x)})$. Our goal is to show that $\theta_1$, $\theta_2$ are constant, and that $\theta_2=2\theta_1$.\\
We already showed that $I(P,Q)=I(|P|,|Q|)$, hence 
$$\int \Big(|\nabla P|_{\gamma_1}^2+|\nabla Q|_{\gamma_1}^2\Big)=\int \Big((\nabla |P|)_{\gamma_1}^2+(\nabla |Q|)_{\gamma_3}^2\Big)$$
and
$$\int \Big(|P|^2|\nabla \theta_1|_{\gamma_1}^2+|Q|^2|\nabla \theta_2|_{\gamma_2}^2\Big)=0.$$
To conclude that $\theta_1$ and $\theta_2$ are constant, one only needs to show that $|P|$ and $|Q|$ do not vanish. To show that $Q$ does not vanish, we use the (real) equation
$$-(4\omega+\beta) Q+\Delta_{\gamma_2} Q=-\frac 12 P^2.$$ 
Noticing that for $L=\Delta_{\gamma_2}-(4\omega+\beta)$, $LQ\leq 0$, we can conclude by using 
the Maximum Principle stated in Theorem 3.5 of \cite{Trudinger}). We can also show that $P$ does not vanish by applying a similar argument to equation 
$$-\omega_{\gamma_1} P+\Delta P=-\frac 12 PQ$$
in a neighborhood of its solution $|P|$.\\
\\
Finally, the relation $\theta_2=2\theta_1$ simply comes from the fact that $$\displaystyle \int |\overline{P}|^2|Q|=\int \overline{P}^2Q,$$ as shown in the proof of (i). \hfill$\blacksquare$

\bigskip

\noindent
{\bf Proof of Theorem \ref{inst}:}\\
For convenience of the notations, we will take $\gamma_1=\gamma_2=1$, although the exact same proof remains valid for arbitrary $\gamma_1,\gamma_2>0$.
In view of Proposition \eqref{positivas}, we may assume that $P,Q>0$. Let
$$\mathcal{L}=\{(u,v)\in H\,:\,M(u,v)=M(P,Q)\}.$$
Following \cite{GR},  it is sufficient to prove the existence of $\Psi\in H$ such that
\begin{enumerate}
\item $\Psi$ is tangent to $\mathcal{L}$ at $(P,Q)$;
\item $J^{-1}\Psi$ is $L^2$-orthogonal to $\partial_{\theta}f(0,0)[P,Q]=i(P,2Q)$ and to\\ $\nabla_{y}f(0,0)[P,Q]=(\nabla P,\nabla Q)$;
\item $\partial_{\theta}f(0,0)[P,Q]$ and $\nabla_{y}f(0,0)[P,Q]$ are linearly independent;
\item  $\langle S''(P,Q)\Psi,\Psi\rangle<0$.
\end{enumerate}
In order for the present paper be self-contained, we briefly explain in the next two steps how these four points can be used to prove Theorem \ref{inst}. For details, we refer the reader to \cite{GR}.

\bigskip

\bigskip

\noindent{\bf Step 1: Construction of an Auxiliary Dynamical System}

\bigskip

\noindent
From conditions 2. and 3., and for some $\epsilon>0$, we build an Auxiliary Dynamical System
$$\mathcal{H}\,:\,\Sigma_{\epsilon}\to \er$$ 
with the following properties:
\begin{itemize}
\item $\forall (U,V)\in \Sigma_{\epsilon}$, $\forall\theta, y\in \er\times \er^n$, $\mathcal{H}(f(\theta,y)v)=\mathcal{H}(v)$;
\item $\forall (U,V)\in \Sigma_{\epsilon}$, $\mathcal{H}'(v)\in H$ and $\mathcal{H'}\,:\,\Sigma_{\epsilon}\to H$ is $C^1$ with bounded derivative;
\item $J\mathcal{H}'(P,Q)=\Psi$.
\end{itemize}
Indeed, consider the mapping
\begin{displaymath}
\begin{array}{lllll}
F\,:\,&H\times \er\times \er^n&\to&\er\\
&((U,V),\theta,y)&\to& \frac 12\|f(\theta,y)(U,V)-(P,Q)\|_2^2.
\end{array}
\end{displaymath}
Using the fact that $\partial_{\theta}f(0,0)[P,Q]$ and $\nabla_{y}f(0,0)[P,Q]$ are linearly independent, one can show, applying the Implicit Theorem Function to $F$ and arguing by convexity, that for $(U,V)$ in a neighbourhood $\mathcal{V}$ of $(P,Q)$, there exists a function $G(U,V)=(G_1(U,V),G_2(U,V))=(\theta(U,V),y(U,V))$ that minimizes $F((U,V),\cdot,\cdot)$ in a ball centered in $(0,0)$; that is, locally, the $L^2$-distance between $(U,V)$ and the orbit of $(P,Q)$ is achieved. Furthermore, one can show that, forall $(\theta,y)\in \er\times\er^n$,
\begin{equation}
\label{G1}
G_1(f(\theta,y)(U,V))\equiv G_1(U,V)-\theta\mod 2\pi
\end{equation}
and
\begin{equation}
\label{G2}
G_2(f(\theta,y)(U,V))=G_2(U,V)-y
\end{equation}
provided that $f(\theta,y)(U,V)\in \mathcal{V}$.
These properties allow to coherently extend the functional 
\begin{displaymath}
\mathcal{H}(U,V)=\langle J^{-1}\Psi,f(G(U,V))(U,V)\rangle
\end{displaymath}
from $\mathcal{V}$ to a entire neighbourhood $\Sigma_{\epsilon}$ of the orbit of $(P,Q)$. Furthermore, in view of \eqref{G1} and \eqref{G2}, it is straightforward that $\mathcal{H}$ is invariant by the action of $f$.\\
Using again the Implicit Function Theorem and the expressions it provides for $G_1'$ and $G_2'$, we can check that $\mathcal{H}'(U,V)\in H$ and that $\mathcal{H}'$ is $C^1$ with bounded derivative.\\  
Finally, from the orthogonality relations expressed in condition 2., one can deduce that $\mathcal{H}'(P,Q)=J^{-1}\Psi$, that is, $J\mathcal{H}'(P,Q)=\Psi$. 

\bigskip

\noindent
{\bf Step 2: Instability}

\bigskip

\noindent
The main idea of the proof is to follow the evolution of the action $S$ along the integral curves of the Auxiliary Dynamical System. More precisely, given $X_0=(U_0,V_0)$ in a neighbourhood $\Sigma_{\epsilon}$ of $\Sigma$ and for a $\sigma>0$, we consider the path
$$\phi\,:\,s\in ]-\sigma,\sigma[\to \phi(X_0,s)\in \Sigma_{\epsilon}$$
such that $\displaystyle \frac d{ds} \phi(X_0,s)=J\mathcal{H}'(X_0,s)$.

\bigskip
\noindent

We consider the evolution of the action along this path, $S(\phi(X_0,s))$. A simple computation then yields
\begin{equation}
\label{derivadas}
\frac d{ds}S(\phi(X_0,s))=P((\phi(X_0,s))\textrm{ and }\frac {d^2}{ds^2}S(\phi(X_0,s))=R((\phi(X_0,s)),
\end{equation}
where
$$P(U,V)=\langle S'(U,V),J\mathcal{H}'(U,V)\rangle$$
and
$$R(U,V)=\langle S''(U,V)i\mathcal{H}'(U,V),J\mathcal{H}'(U,V)\rangle+\langle S'(U,V),J\mathcal{H}''(U,V)(J\mathcal{H}'(U,V))\rangle.$$
Using the Taylor expansion, we obtain the existence of $\xi\in[0,1]$ such that
\begin{equation}
\label{Taylor}
S(\phi(X_0,s))=S(X_0)+P(X_0)s+\frac 12 R(\phi(X_0,\xi s))s^2.
\end{equation}
Noticing that $S'(P,Q)=0$ (from \eqref{elliptic}) and $J\mathcal{H}'(P,Q)=\Psi$, we obtain that
$$R(P,Q)=\langle S''(P,Q)\Psi,\Psi\rangle<0,$$
yet, from \eqref{Taylor}, for $X_0$ in a neighbourdhood of $(P,Q)$ and for small $s$, 
\begin{equation}
\label{Taylor2}
S(\phi(X_0,s))\leq S(X_0)+P(X_0)s.
\end{equation}
By intersecting the manifold $\mathcal{W}_{(P,Q)}$ with the trajectories of the Auxiliary Dynamical System, using the Implicit Theorem Function, it is possible to obtain a uniform version of \eqref{Taylor2}, namely, for some $\epsilon>0$,
\begin{equation}
\label{Taylor3}
\forall X_0\in\Sigma_{\epsilon},\,\exists s\in]-\sigma,\sigma[,\,S(\phi(X_0,s))\leq S(X_0)+P(X_0)s.
\end{equation}
This means that $P$ measures the variations of $S$ (hence of $E$) along the trajectories pf the Auxiliary Dynamical System.\\
The crucial step is now to prove that $P$ also measures the variations of $\mathcal{H}$ along the the flow of the initial system \eqref{quadratic}. More precisely, considering the solution $(u(t),v(t))$ of \eqref{quadratic} with initial data $X_0$, we have
\begin{equation}
\label{trajin}
\frac d{dt} \mathcal{H}(u(t),v(t))=-P(u(t),v(t)).
\end{equation}

\bigskip

\noindent
This can be achieved by justifying the following formal computation:
$$\mathcal{H}(u(t),v(t))-\mathcal{H}(u_0,v_0)=\int_0^t \langle \mathcal{H}'(u(\tau),v(\tau)),(u_t(\tau),v_t(\tau)) \rangle d\tau$$
$$=\int_0^t\langle  \mathcal{H}'(u(\tau),v(\tau)),J^{-1} E'(u(t),v(t))\rangle d\tau=-\int_0^t P(u(\tau),v(\tau))d\tau$$
since mass is conserved along the trajectories of the Auxiliary Dynamical System \eqref{quadratic} and $J^{-1}$ is skew-adjoint.

\bigskip

\noindent
Finally, setting 
$$\mathcal{P}=\{(U,V)\in \Sigma_{\epsilon}\,:\,S(U,V)<S(P,Q)\textrm{ and }P(U,V)\neq 0\},$$
it can be shown that, for $X_0\in\mathcal{P}$, $P(u(t),v(t))$ remains bounded away of the origin as long as the solution $(u(t),v(t))$ exists. This implies that solutions of \eqref{quadratic} for initial data $X_0\in\mathcal{P}$ must leave in finite time any neighbourhood of $\Sigma$.
Indeed, 
$$\Big|\frac d{dt}\mathcal{H}(u(t),v(t))\Big|=|P(u(t),v(t)|\geq C(X_0)>0,$$
which contradicts the fact that  $\mathcal{H}$ in bounded in any neighbourhood of $\Sigma$. (recall $\mathcal{H}$ is invariant by $f(\theta,y)$, $\theta\in\er$ $y\in \er^n$). 

\bigskip

\noindent
Now, following the action along the trajectory of the Auxiliary Dynamical System that contains $(P,Q)$ it can be shown that $\mathcal{P}$ contains points arbitrarely close to $(P,Q)$ of the form $\phi((P,Q),s)$, that is, belonging to the considered trajectory.

\bigskip

\noindent
Also, setting $\displaystyle W(U,V)=Re\int \overline{U}^2V$ the potential energy and observing that the map
$$A\,:\,s\to W(\phi((P,Q),s))$$
is $C^1$ and has a non vanishing derivative at the origin, for small $s$ with the adequate sign,
$$W(\phi((P,Q),s))<W(P,Q).$$
Putting $X_0=\phi((P,Q),s)$ and considering the solution $X(t)=(u(t),v(t))$ of \eqref{quadratic} with initial data $X_0$, we have, as long as the solution exists, $$W(X(t))<W(P,Q).$$ Indeed, if at some point $W(X(t))=W(P,Q)$ then we would obtain a contradiction with the fact that $X_0\in\mathcal{P}$ and that $(P,Q)$ is a solution of \eqref{minimization}:
$$X(t)\in \mathcal{W}_{(P,Q)}\textrm{ and }S(X(t))\leq S(X_0)<S(P,Q).$$
\noindent
Since $E$ is conserved by the flow of \eqref{quadratic}, this is enough to prove that $(u(t),v(t))$ is bounded in $H$ and global.\hfill$\blacksquare$

\bigskip

\bigskip

\noindent
{\bf End of the Proof of Theorem \ref{inst}:}\\

\noindent
We now exhibit $\Psi\in H$ satisfying properties $1.$, $2.$, $3.$ and $4.$\\
Let $\e>0$. We begin by considering the curve  
$$\Gamma\,:\, t\in[0,\e[\to \Big(\gamma(t)\lambda^{\frac n2}(t)P(\lambda(t)\cdot) , \alpha(t)\lambda^{\frac n2}(t)Q(\lambda(t)\cdot)\Big)$$
where $\alpha, \gamma$ and $\lambda$ are smooth real functions to be chosen later, such that
\begin{equation}
\label{cond1}
\alpha(0)=\gamma(0)=\lambda(0)=1\,\quad(\textrm{ that is, }\Gamma(0)=(P,Q)).
\end{equation}

\bigskip

\noindent
1. Setting
$$k=\frac{\int P^2}{4\int Q^2},$$
the condition
\begin{equation}
\label{cond2}
\gamma^2k+\alpha^2=k+1
\end{equation}
assures that $\Gamma\subset\mathcal{L}$, and, in particular,
\begin{equation}
\Psi=\Gamma'(0)
\end{equation}
is tangent to $\mathcal{L}$ at $(P,Q)$. 

\bigskip

\noindent
2. Noticing that
$$\Psi=((\lambda^{\frac n2}\gamma)'(0)P+\lambda'(0)\nabla P,(\lambda^{\frac n2}\alpha)'(0)Q+\lambda'(0)\nabla Q)$$ has real components,
$$i\Psi \perp \nabla_{y}f(0,0)[P,Q].$$
Also,  $i\Psi\perp i(P,2Q)$ since $\Psi\in T_{\mathcal{L}}(P,Q)$ and $(P,2Q)\perp T_{\mathcal{L}}(P,Q)$.

\bigskip

\noindent
3. Since $i(P,2Q) \perp (\nabla P, \nabla Q)$, these two vectors are linearly independent. 

\bigskip

\noindent
4. We begin by computing the energy (\ref{energy}) along the path $\Gamma$:
\begin{displaymath}
\begin{array}{llll}
E(\Gamma(t))&=&\displaystyle\gamma^2\lambda^2\int |\nabla P|^2+\alpha^2\lambda^2\int |\nabla Q|^2+\beta\alpha^2\int Q^2-\gamma^2\alpha\lambda^{\frac n2}\int P^2Q\\
\\
&=&\displaystyle\frac 1k\left((k+1-\alpha^2)\lambda^2\int |\nabla P|^2+k\alpha^2\lambda^2\int |\nabla Q|^2+k\beta\alpha^2\int Q^2\right.\\
\\
&&\displaystyle -(k+1-\alpha^2)\alpha\lambda^{\frac n2}\int P^2Q.
\end{array}
\end{displaymath}
Differentiating with respect ro $t$,
\begin{displaymath}
\begin{array}{llll}
k\displaystyle\frac{d}{dt}E(\Gamma(t))&=&\alpha'A(t)+\lambda'B(t),
\end{array}
\end{displaymath}
with
$$A(t)=-2\alpha\lambda^2\int |\nabla P|^2+2k\alpha\lambda^2\int |\nabla Q|^2+2k\beta\alpha\int Q^2+(3\alpha^2-k-1)\lambda^{\frac n2}\int P^2Q$$
and
$$B(t)=2\lambda(k+1-\alpha^2)\int |\nabla P|^2+2k\alpha^2\lambda\int |\nabla Q|^2+\frac n2(\alpha^2-k-1)\alpha\lambda^{\frac {n-2}2}\int P^2Q.$$
Now, observe that since $(P,Q)$ is a solution of \eqref{elliptic}, $A(0)=B(0)=0$\\
(see \cite{Saut}, (5.2)). Hence, putting $\alpha_0=\alpha'(0)$ and $\lambda_0=\lambda'(0)$,
\begin{displaymath}
\begin{array}{llll}
k\displaystyle\frac{d^2}{dt^2}E(\Gamma(t))|_{t=0}&=&\alpha_0A'(0)+\lambda_0B'(0)\\
\\
&=&\displaystyle \alpha_0^2\Big(-2\int|\nabla P|^2+2k\int|\nabla Q|^2+2k\beta\int Q^2+6\int P^2Q\Big)\\
\\
&+&\displaystyle \lambda_0^2\Big(2k\int|\nabla P|^2+2k\int|\nabla Q|^2-\frac{n(n-2)}{4}k\int P^2Q\Big)\\
\\
&+&\displaystyle 2\alpha_0\lambda_0\Big(-4\int|\nabla P|^2+4k\int|\nabla Q|^2 -\frac{(k-2)}{2}n\int P^2Q\Big).
\end{array}
\end{displaymath}
Using once again that $A(0)=B(0)=0$, this quantity can be re-written in terms of $\displaystyle \int P^2Q$ and $\displaystyle \int Q^2$ exclusively:
\begin{displaymath}
\begin{array}{llll}
k\displaystyle\frac{d^2}{dt^2}E(\Gamma(t))|_{t=0}&=&\displaystyle \alpha_0^2\Big((k+4)\int P^2Q\Big)+\lambda_0^2\Big(\frac{n(4-n)}{4}k\int P^2Q\Big)\\
\\
&+&\displaystyle 2\alpha_0\lambda_0\Big(-4k\beta\int Q^2 +\frac{(k-2)(4-n)}{2}\int P^2Q\Big).
\end{array}
\end{displaymath}

\bigskip

\noindent
The determinant of $\displaystyle\frac{d^2}{dt^2}E(\Gamma(t)|_{t=0}$ as a quadratic form in $(\alpha_0,\lambda_0)$ is given by
$$\Delta=(k+4)\frac{n(4-n)}{4}\Big(\int P^2Q\Big)^2-\Big(\frac{(k-2)(4-n)}{2}\int P^2Q-4k\beta\int Q^2\Big)^2.$$
For $d\geq 4\,(n\geq 5)$, $\Delta<0$.
For $d=3\,(n=4)$ and $\beta\neq 0$,
$$\Delta=-16\beta^2k\int Q^2<0.$$
Hence, in both these situations, one can choose $\alpha(t),\lambda(t)$ such that $$\displaystyle\frac{d^2}{dt^2}E(\Gamma(t))|_{t=0}<0.$$

\noindent
Now, observe that
$$\frac{d}{dt}S(\Gamma(t))=\langle S'(\Gamma(t)),\Gamma'(t)\rangle$$
and
$$\frac{d^2}{dt^2}S(\Gamma(t))=\langle S'(\Gamma(t)),\Gamma''(t)\rangle+\Gamma'(t)^T\Big[S''(\Gamma(t))\Big]\Gamma'(t).$$
Since $S'(\Gamma(0))=S'(P,Q)=0$, setting $t=0$ yields
$$\frac{d^2}{dt^2}S(\Gamma(t))|_{t=0}=\langle S''(P,Q)\Psi,\Psi\rangle.$$
Finally $\Gamma\subset\mathcal{L}$, $\displaystyle \frac{d}{dt}M(\Gamma(t))=0$ and $$\langle S''(P,Q)\Psi,\Psi\rangle=\frac{d^2}{dt^2}E(\Gamma(t))|_{t=0}<0.$$\hfill$\blacksquare$

\section{A stability result}\label{S4}
\label{stability}

Define, for any $(u,v)\in H$ such that $J(u,v)>0$,
$$
GN(u,v)=\frac{M(u,v)^{\frac{3}{2}-\frac{d+1}{4}}K(u,v)^{\frac{d+1}{4}}}{J(u,v)}.
$$
This functional is closely related with a vector-valued Gagliardo-Nirenberg inequality: if one sets
\begin{equation}\label{min_gagliardo}
C_{GN}^{-1}:=\inf\{GN(u,v): J(u,v)>0\},
\end{equation}
then $C_{GN}$ is the optimal constant of the inequality
$$
Re \int \overline{u}^2 v \le C\left(\int |u|^2 + 4|v|^2\right)^{\frac{3}{2}-\frac{d+1}{4}}\left(\int |\nabla u|^2 + |\nabla v|^2\right)^{\frac{d+1}{4}}.
$$
\begin{Lem}\label{lemma_GN}
Suppose that $\beta=0$. Then the set of solutions for the minimization problem \eqref{min_gagliardo} is $G$, up to scalar multiplication and scaling.
\end{Lem}
\noindent
\textbf{Proof:}\\
By \cite{Saut}, we know that $G\neq\emptyset$ is the set of solutions of \eqref{minimization}. Let $\mathcal{Q}\in G$ and $W=(w,z)$ be such that $J(W)>0$. Recall that $I(\mathcal{Q})=J(\mathcal{Q})>0$. Define
$$
\nu = \left(\frac{J(\mathcal{Q})M(W)}{M(\mathcal{Q})J(W)}\right)^{\frac{1}{2p}}
$$
and
$$
\zeta = \left(\nu^2\left(\frac{M(W)}{M(\mathcal{Q})}\right)\right)^{\frac{1}{N}}.
$$
Then $Z(x)=\nu W(\zeta x)$ satisfies
$$
J(Z)=J(\mathcal{Q}), \quad M(Z)=M(\mathcal{Q}) \quad GN(Z)=GN(W).
$$
By the minimality of $\mathcal{Q}$, $I(\mathcal{Q})\le I(Z)$, which implies that $GN(\mathcal{Q})\le GN(Z)=GN(W)$. Therefore $\mathcal{Q}$ is a solution of  \eqref{min_gagliardo}. On the other hand, if $W$ is a solution of \eqref{min_gagliardo}, then one has necessarily $GN(Z)=GN(\mathcal{Q})$, which implies that $I(Z)=I(\mathcal{Q})$. Therefore $Z\in G$, which concludes our proof.\hfill$\blacksquare$

\begin{Lem}\label{equivalenciasubcritico}
Suppose that $d\le 2$ and $\beta=0$. Then the set of ground states $G$ is the set of solutions of \eqref{min_energia}.
\end{Lem}
\noindent
\textbf{Proof:}\\
Let $W=(w,z)\in H$ be such that $M(W)=\nu$. For any $\lambda>0$, define $W_\lambda(x):=\lambda^{\frac{d+1}{2}}W(\lambda x)$. Consider the function
$$
\lambda\mapsto f(\lambda)=E(W_\lambda),\ \lambda>0
$$
Since $d\le 2$, $f$ has a unique minimum $\lambda_0$. Let $Z=W_{\lambda_0}$. Then $f'(\lambda_0)=0$, which implies that 
$$
K(Z)=\frac{d+1}{6}J(Z).
$$
Therefore,
\begin{equation}\label{energia_z}
E(Z)=\frac{d-3}{2d+2}K(Z).
\end{equation}

Let $\mathcal{Q}\in G$. Notice that, by Pohozaev's equality, the same relations are valid for $\mathcal{Q}$:
\begin{equation}\label{energia_q}
K(\mathcal{Q})=\frac{d+1}{6}J(\mathcal{Q}),\quad E(\mathcal{Q})=\frac{d-3}{2d+2}K(\mathcal{Q}).
\end{equation}
By Lemma \ref{lemma_GN}, we have
$$
\frac{(d+1)\nu^{\frac{3}{2}-\frac{d+1}{4}}K(\mathcal{Q})^{\frac{d+1}{4}}}{6K(\mathcal{Q})}=GN(\mathcal{Q})\le GN(Z) = \frac{(d+1)\nu^{\frac{3}{2}-\frac{d+1}{4}}K(Z)^{\frac{d+1}{4}}}{6K(Z)}
$$
and so $K(Z)\le K(\mathcal{Q})$. Hence, by \eqref{energia_z} and \eqref{energia_q},
$$
E(W)\ge E(Z) =\frac{d-3}{2d+2}K(Z) \ge \frac{d-3}{2d+2}K(\mathcal{Q}) = E(\mathcal{Q}),
$$
and so $\mathcal{Q}$ is a solution of \eqref{min_energia}.

On the other hand, if $W$ is also a solution of \eqref{min_energia}, then one must have $Z=W$ and $K(W)=K(\mathcal{Q})$. Again by \eqref{energia_z} and \eqref{energia_q}, 
$$
J(W)=\frac{6}{d+1}K(W)=\frac{6}{d+1}K(\mathcal{Q})=J(\mathcal{Q}).
$$
Moreover, since $M(W)=M(\mathcal{Q})$, one has $I(W)=I(\mathcal{Q})$. This implies that $W$ is a solution of \eqref{minimization}, \textit{i.e.}, $W\in G$.\hfill$\blacksquare$

\textbf{Sketch of the proof of Theorem \ref{Theo:St}:}
The proof follows the same steps as \cite[Proposition 4]{Saut}: suppose, by contradiction, that there exist sequences $\{(u^n_0,v^n_0)\}_{n\in\mathbb{N}}\subset H$ and $\{t_n\}_{n\in\mathbb{N}}\subset\R^+$ with
\begin{equation}\label{dados_convergem}
\inf_{(P,Q)\in G}\|(u_0^n,v_0^n)-(P,Q)\|_{H}\to 0
\end{equation}
and such that the corresponding solutions $(u^n,v^n)$ satisfy
\begin{equation}\label{solucoes_nao_convergem}
\inf_{(P,Q)\in G}\|(u^n(t_n),v^n(t_n))-(P,Q)\|_{H}>\delta.
\end{equation}
By \eqref{dados_convergem}, for any given $(P,Q)\in G$,
$$
M(u_0^n,v_0^n)\to M(P,Q),\quad E(u_0^n,v_0^n)\to E(P,Q).
$$
Using the conservation of mass and energy, one has
$$
M(u^n(t_n),v^n(t_n))\to M(P,Q),\quad E(u^n(t_n),v^n(t_n))\to E(P,Q).
$$
This implies that (up to a normalization) $\{(u^n(t_n),v^n(t_n))\}_{n\in\mathbb{N}}$ is a minimizing sequence of problem \eqref{min_energia}. The argument of \cite{Saut} implies that 
$$(u^n(t_n),v^n(t_n))\to (\tilde{P}, \tilde{Q}),
$$
where $(\tilde{P}, \tilde{Q})$ is a solution of \eqref{min_energia}, that is, $(\tilde{P}, \tilde{Q})\in G$. This convergence contradicts \eqref{solucoes_nao_convergem}, thus finishing the proof.\hfill$\blacksquare$


\bigskip

\noindent
{\bf Acknowledgements}
This work was developed in the frame of the CAPES-FCT convenium {\it Equa\c{c}\~oes de evolu\c{c}\~ao dispersivas}. 
Sim\~ao Correia and Filipe Oliveira would like to thank the kind hospitality of IMPA, Instituto de Matem\'atica Pura e Aplicada, and of the Institute of Mathematics at the Federal University of Rio de Janeiro. Ad\'an Corcho would like to thank the kind hospitality of Instituto Superior T\'ecnico.
Sim\~ao Correia was partially supported by Funda\c{c}\~ao para a Ci\^encia e Tecnologia, through the grant SFRH/BD/96399/2013 and through contract
UID/MAT/04561/2013.\\ Filipe Oliveira  was partially supported by the Project CEMAPRE - UID/ MULTI/00491/2013 financed by FCT/MCTES through national funds.

\end{document}